\documentclass[12pt]{article}

\usepackage{amssymb}
\usepackage{graphics} 
\usepackage{url}
\usepackage{hyperref}
\usepackage[authoryear,round]{natbib}
\usepackage{amssymb,color}
\usepackage{graphics,amsmath}

\begin{document}

\newenvironment {proof}{{\noindent\bf Proof.}}{\hfill $\Box$ \medskip}

\newtheorem{theorem}{Theorem}[section]
\newtheorem{lemma}[theorem]{Lemma}
\newtheorem{condition}[theorem]{Condition}
\newtheorem{proposition}[theorem]{Proposition}
\newtheorem{remark}[theorem]{Remark}
\newtheorem{hypothesis}[theorem]{Hypothesis}
\newtheorem{corollary}[theorem]{Corollary}
\newtheorem{example}[theorem]{Example}
\newtheorem{definition}[theorem]{Definition}

\renewcommand {\theequation}{\arabic{section}.\arabic{equation}}

\def \non {{\nonumber}}
\def \noin {{\noindent}}
\def \hat {\widehat}
\def \tilde {\widetilde}
\def \E {\mathbb{E}}
\def \P {\mathbb{P}}
\def \R {\mathbb{R}}
\def \N {\mathbb{N}}
\def \I {\mathbb{I}}

\title{\large {\bf Existence and uniqueness }
{\bf of  reflecting diffusions in cusps}}
                                                       
\author{\begin{tabular}{ll}                               
Cristina Costantini & Thomas G. Kurtz \\
Dipartimento di Economia & Department of Mathematics\\
 and Unit\`a locale INdAM & and Department of Statistics\\
Universit\`a di Chieti-Pescara & University of Wisconsin - Madison \\            
v.le Pindaro 42 & 480 Lincoln Drive  \\                                                    
65127 Pescara, Italy & Madison, WI  53706-1388, USA \\                    
c.costantini@unich.it & kurtz@math.wisc.edu    \\
& \url{http://www.math.wisc.edu/~kurtz/}  \\                   
\end{tabular}}

\date{October 12, 2017}

\maketitle

\begin{abstract}
We consider stochastic differential equations with 
(oblique) reflection in a $2$-dimensional domain that has a 
cusp at the origin, i..e. in a neighborhood of the origin 
has the form $\{(x_1,x_2):0<x_1\leq\delta_0,\psi_1(x_1)<x_2<\psi_
2(x_1)\}$, with 
$\psi_1(0)=\psi_2(0)=0$, $\psi_1'(0)=\psi_2'(0)=0$. 

Given a vector field $\gamma$ of directions of reflection 
at the boundary points other than the origin, 
defining directions of reflection at the origin 
$\gamma^i(0):=\lim_{x_1\rightarrow 0^{+}}\gamma (x_1,\psi_i(x_1))$, $
i=1,2,$ and assuming there 
exists a vector $e^{*}$ such that $\langle e^{*},\gamma^i(0)\rangle 
>0$, $i=1,2$, and 
$e^{*}_1>0$, we prove weak existence 
and uniqueness of the solution starting at the origin and 
strong existence and uniqueness starting away from the 
origin. 

Our proof uses a new scaling result and a coupling argument.
\vspace{.1in}

\vspace{.1in}


\end{abstract}

\setcounter{equation}{0}
\section{Introduction}

In this work we prove existence and uniqueness of 
reflecting diffusions in $2$-dimensional domains 
with cusps.  By saying that the domain has a cusp, we 
mean that in a neighborhood of some point, which we 
take to be the origin, the domain has the form 
\[\{(x_1,x_2)\in D:0<x_1\leq\delta_0\}=\{(x_1,x_2):0<x_1\leq\delta_
0,\psi_1(x_1)<x_2<\psi_2(x_1)\},\]
where we assume $\psi_1$ and $\psi_2$ are $C^1$ with
\[\psi_1(0)=\psi_2(0)=0,\quad\psi_1'(0)=\psi_2'(0)=0,\]
and, in general, the boundary is $C^1$ away from the origin. 

The direction of 
reflection, $\gamma$, is assumed Lipschitz continuous 
on the smooth part of the boundary, with a uniformly 
positive scalar product with the inward normal.  At the 
tip, 
\[\gamma^i(0):=\lim_{x_1\rightarrow 0^{+}}\gamma (x_1,\psi_i(x_1)
),\quad i=1,2,\]
is assumed to exist, and for some  $e^{*}\in {\Bbb R}^2$ 
\[\langle e^{*},\gamma\rangle >0,\quad\forall\gamma\in\left\{\left
[\begin{array}{c}
1\\
0\end{array}
\right],\gamma^1(0),\gamma^2(0)\right\}.\]
See Section \ref{assumptions} for the complete formulation of 
our assumptions. 

In the case of a domain $D$ of the form 
\[D:=\{(x_1,x_2):\,0<x_1,\,\,\psi_1(x_1):=-x_1^{\beta_1}<x_2<\psi_
2(x_1):=x_1^{\beta_2}\},\]
(where either $\beta_1=\beta_2>1$ or $\beta_1>2\beta_2-1$, $\beta_
2>1$), under the 
assumption that on each of $\{(x_1,x_2):\,0<x_1,\,x_2=-x_1^{\beta_
1}\}$ and 
$\{(x_1,x_2):\,0<x_1,\,x_2=x_1^{\beta_2}\}$ the direction of reflection forms 
a constant angle with the inward normal, weak existence 
and uniqueness of reflecting Brownian motions 
have been exhaustively studied by \cite{DBT93a}.  
\cite{DBT93b} gives a complete characterization of the 
cases in which, in the above setup, the reflecting 
Brownian motion is a semimartingale.  For general 
continuous $\psi_1$, $\psi_2$ $(\psi_1(0)=\psi_2(0)=0$, $\psi_2(x_
1)>\psi_1(x_1)$ for 
every $x_1>0$), the case when the directions of reflection 
on $\{(x_1,x_2):\,0<x_1,\,x_2=\psi_1(x_1)\}$ and 
$\{(x_1,x_2):\,0<x_1,\,x_2=\psi_2(x_1)\}$ are constant, opposite vertical 
vectors - a case when the process is not a 
semimartingale - has been studied by \cite{BT95} and 
\cite{BKR09}.  In higher dimensions,
normally reflecting diffusions in domains 
with H\"older cusps have 
been studied by \cite{FT96} by analytical techniques.  

Here we characterize the reflecting 
diffusion as the solution of a stochastic differential 
equation with reflection (SDER) which will always 
be a semimartingale. In particular, we recover the results 
by \cite{DBT93a} and \cite{DBT93b} for the cases when 
the process is a semimartingale, except for the case 
when $\gamma^1(0)$ and $\gamma^2(0)$ point at each other and $\beta_
2<2$.

First, we show that our conditions imply
that, starting away 
from the origin, the origin is never reached. Therefore 
we easily obtain strong existence and uniqueness of the 
reflecting diffusion from known 
results on existence and uniqueness in smooth domains 
(Section \ref{off0}). 

Moreover, the fact that, starting away from the origin, 
the reflecting diffusion is well defined for all times 
allows us to obtain a weak solution of the SDER starting 
at the origin as the limit of solutions starting away 
from the origin (Section \ref{tipexist}).  To this end, we 
employ a random time change of the SDER (the same 
that is used in \cite{Kur90} to obtain a solution of a 
patchwork martingale problem from a solution of the 
corresponding constrained martingale problem) that 
makes it particularly simple to prove relative 
compactness of the processes.  

The main result of this paper, however, is weak 
uniqueness of the solution to the SDER starting at the 
origin (Section \ref{tipunq}).  Our assumptions on the 
direction of reflection guarantee that any solution 
starting at the origin immediately leaves it.  Since the 
distribution of a solution starting away from the origin 
is uniquely determined, the distribution of a solution 
starting at the origin is determined by its exit 
distribution from an arbitrarily small neighborhood of 
the origin.  The crucial ingredient that allows us to 
understand the behavior of the process near the origin is 
a scaling result (Section \ref{tipscale}).  Combined with a 
coupling argument based on \cite{LR86}, this scaling 
result shows that indeed all solutions starting at the 
origin must have the same exit distribution from every 
neighborhood of the origin.  For a more detailed 
discussion of our approach, see the beginning of Section 
\ref{0}.  Some technical lemmas that are needed in our 
argument are proved in Section \ref{technic}.  The Feller 
property is proved in Section \ref{Feller}.  

The most general uniqueness result for SDER in 
piecewise ${\cal C}^1$ domains can be found in \cite{DI93}.  
Reflecting diffusions in piecewise smooth 
domains are characterized as solutions of constrained 
martingale problems in \cite{Kur90}, and \cite{CK15} 
reduces the problem of proving uniqueness for the 
solution of a constrained martingale problem (as well as 
of a martingale problem in a general Polish space) to 
that of proving a comparison principle for viscosity 
semisolutions of the corresponding resolvent equation.  
None of these results applies to the situation we are 
considering here.  In particular, \cite{DI93} makes the 
assumption that the convex cone generated by the 
normal vectors at each point does not contain any 
straight line, which is violated at the tip of the cusp.  

Finally, we wish to mention that our work was partly 
motivated by diffusion approximations for some 
queueing models where domains with cusplike 
singularities appear (in particular \cite{KKLW09}). 
These models are in higher dimensions, but this 
paper is intended as a first contribution in the direction of 
understanding reflecting diffusions in such 
domains. 

The notation used in the paper is collected in Section 
\ref{notation}.

\setcounter{equation}{0}

\section{Formulation of the problem and 
assumptions}\label{assumptions}

We are interested in studying diffusion processes 
with oblique reflection in the closure of a simply 
connected  
$2$-dimensional 
domain $D\subset [0,\infty )\times\R$ with a boundary $\partial D$ that is $
C^1$ 
except at a single point 
(which we will take to be the origin $0$), where 
the domain has a cusp. More precisely $D$ satisfies the 
following. 

\begin{condition}\label{domcond} 
\item[a)]$D$ is a bounded, simply connected domain in 
$[0,\infty )\times\R$ with $0\in\partial D$.

\item[b)]$\partial D$ is $C^1$except at $0$.

\item[c)]There exists a $\delta_0>0$ and continuously 
differentiable functions $\psi_1$ and $\psi_2$ with $\psi_1\leq\psi_
2$ and 
\[\psi_1(0)=\psi_2(0)=0,\quad\psi_1'(0)=\psi_2'(0)=0\]
such that 
\[\{(x_1,x_2)\in D:x_1\leq\delta_0\}=\{(x_1,x_2):0<x_1\leq\delta_
0,\psi_1(x_1)<x_2<\psi_2(x_1)\},\]
and
\[\lim_{x_1\rightarrow 0+}\frac {\psi_1(x_1)}{\psi_2(x_1)-\psi_1(
x_1)}=L\in (-\infty ,\infty ).\]
\end{condition}

The direction of reflection is assigned at all points of 
the boundary except the origin and is given by a 
unit vector field $\gamma$ verifying the following condition. 

\begin{condition}\label{dircond}
\item[a)]$\gamma$$:\partial D-\{0\}\rightarrow {\Bbb R}^2$ is locally Lipschitz continuous and satisfies 
\[\inf_{x\in\partial D-0}\langle\gamma (x),\nu (x)\rangle >0.\]
$ $The mappings 
\[x_1\in (0,\delta_0]\rightarrow\gamma^i(x_1):=\gamma (x_1,\psi_i
(x_1)),\quad i=1,2,\]
are Lipschitz continuous and hence the limits 
\[\gamma^i(0):=\lim_{x_1\rightarrow 0^{+}}\gamma (x_1,\psi_i(x_1)
),\quad i=1,2,\]
exist. 

\item[b)]Let $\Gamma (0)$ be the convex cone generated by $ $
\[\left\{\left[\begin{array}{c}
1\\
0\end{array}
\right],\gamma^1(0),\gamma^2(0)\right\}.\]
There exists $e^{*}\in {\Bbb R}^2$ such that 
\[\langle e^{*},\gamma\rangle >0,\quad\forall\gamma\in\Gamma (0).\]
Of course, without loss of generality, we can suppose 
that $|e^{*}|=1$.
\end{condition}

\begin{remark}
Condition \ref{dircond}(b) can be reformulated as follows. 
In a neighborhood of the origin, we can 
view $D$ as being the intersection of three ${\cal C}^1$ domains, 
\[\{x:\,x_2>\psi_1(x_1)\},\quad \{x:\,x_2<\psi_2(x_1)\},\quad \{x
:\,x_1>0\},\]
with unit inward normal vector at the origin, 
respectively, 
\[\nu^1(0)=\left[\begin{array}{c}
0\\
1\end{array}
\right],\quad\nu^2(0)=\left[\begin{array}{c}
0\\
-1\end{array}
\right],\quad\nu^0(0)=\left[\begin{array}{c}
1\\
0\end{array}
\right].\]
Then, letting the normal cone at the origin, $N(0)$, be the 
closed, convex cone generated by $\{\nu^1(0),\nu^2(0),\nu^3(0)\}$, 
Condition \ref{dircond}(b) 
is equivalent to requiring that there exists $e^{*}\in N(0)$ 
such that 
\[\langle e^{*},\gamma\rangle >0,\quad\forall\gamma\in\Gamma (0),\]
where we can think of $\Gamma (0)$ as the closed, convex cone generated 
by the directions of 
reflection at the origin for each of the three domains. 
In other terms, Condition \ref{dircond}(b) is the analog 
of the condition usually assumed in the literature for 
polyhedral domains (see e.g. \cite{VW85}, \cite{TW93} or 
\cite{DW96}). Note that, in contrast, the condition that 
there exists $e^{**}\in\Gamma (0)$ such that 
\[\langle e^{**},\nu\rangle >0,\quad\forall\nu\in N(0),\]
can never be satisfied at a cusp, because $\nu^2(0)=-\nu^1(0)$. 
\end{remark}

\begin{remark}\label{cone}
Note that, under Condition \ref{dircond}(b), $\left[\begin{array}{c}
1\\
0\end{array}
\right]$ can be 
expressed as a positive linear
combination of $\gamma^1(0)$ and $\gamma^2(0)$, so 
that 
$\Gamma (0)$ coincides with the closed, convex cone generated by $ $
\[\left\{\gamma^1(0),\gamma^2(0)\right\}.\]
\end{remark}

We seek to  characterize the diffusion process with 
directions of reflection $\gamma$ as the solution of a stochastic 
differential equation  driven by a standard Brownian 
motion $W$:
\begin{eqnarray}
X(t)=X(0)+\int_0^tb(X(s))ds&+&\int_0^t\sigma (X(s))dW(s)+\int_0^t
\gamma (s)d\Lambda (s),\qquad t\geq 0,\non\\
\gamma (s)\in\Gamma_1(X(s)),\quad&&\quad d\Lambda -a.e.,\qquad\qquad 
t\geq 0,\label{SDER}\\
X(t)\in\bar {D},\quad&&\int_0^t{\bf 1}_{\partial D}(X(s))d\Lambda 
(s)=\Lambda (t),\quad\qquad t\geq 0,\non\end{eqnarray}
where $\Lambda$ is nondecreasing, $\Gamma_1(0)$ is the convex hull of 
$\gamma^1(0)$ and $\gamma^2(0)$ and for $x\in\partial D-\{0\}$ $\Gamma_
1(x):=\{\gamma (x)\}$ and 
$\gamma$ is almost surely measurable.

We make the following assumptions on the coefficients.

\begin{condition}\label{eqcond}
\item[a)]$\sigma$ and $b$  are Lipschitz continuous on $\overline 
D$.

\item[b)]$(\sigma\sigma^T)(0)$ is nonsingular.
\end{condition}

We will denote 
\begin{equation}Af(x):=Df(x)b(x)+\frac 12\mbox{\rm tr}((\sigma\sigma^
T)(x)D^2f(x)).\label{op}\end{equation}
\begin{definition}
A stochastic process $Z$ is {\em compatible\/} with a Brownian 
motion $W$ if 
for each $t\geq 0$, 
$W(t+\cdot )-W(t)$ is independent of ${\cal F}_t^{W,Z}$, where $\{
{\cal F}^{W,Z}_t\}$ is 
the filtration generated by $W$ and $Z$.
\end{definition}

\begin{definition}
Given a standard Brownian motion $W$ and $X(0)\in\overline D$ 
independent of $W$, $(X,\Lambda )$ is a {\em strong solution\/} of 
(\ref{SDER}) if 
 $(X,\Lambda )$ is 
adapted to the filtration generated by $X(0)$ and $W$
and the equation is satisfied.

$(X,\Lambda ,W)$, defined on some probability space,
 is a {\em weak solution\/} of (\ref{SDER}) if  
$W$ is a standard Brownian motion,
$(X,\Lambda )$ is compatible 
with $W$,  and the equation 
is satisfied.  

Given an initial distribution $\mu\in {\cal P}(\overline D)$, {\em weak }
{\em uniqueness\/} or {\em uniqueness in distribution\/} holds if for all weak solutions with 
$P\{X(0)\in\cdot \}=\mu$,  $X$ has the same distribution on ${\cal C}_{\overline 
D}[0,\infty )$. 

{\em Strong uniqueness\/} holds if for any standard Brownian 
motion $W$ and weak solutions 
$(X,\Lambda ,W)$, $(\tilde {X},\tilde{\Lambda },W)$ such that $X(
0)=\tilde {X}(0)$ a.s. and
  $(X,\Lambda ,\tilde {X},\tilde{\Lambda })$ is 
compatible with $W$, $X=\tilde {X}$ a.s.
\end{definition}

\begin{remark}
Of course, any strong solution is a weak solution. 
Existence of a weak solution and strong uniqueness 
imply that the weak solution is a strong solution (c.f., 
\cite{YW71} and \cite{Kur14}).
\end{remark}

For processes starting 
away from the tip, existence and uniqueness follows 
from results of \cite{DI93} and the fact that under 
our conditions, the solution never hits the tip.  For 
processes starting at the tip, we only prove weak 
existence and uniqueness.  The proof is based on 
rescaling of the process near the tip and a coupling 
argument.

\setcounter{equation}{0}

\section{Strong existence and uniqueness starting at 
$x^0\neq 0$}\label{off0}

Our first result is that, for every 
$x^0\in\overline D-\{0\}$, $(\ref{SDER})$ has a unique  
strong solution with $X(0)=x^0$, well-defined for all times. 
In fact, by 
\cite{DI93}, for each $n>0$, the solution, $X$, is well-defined up to 
\begin{equation}\tau_n:=\inf\{t\geq 0:\,X_1(t)<\frac 1n\},\label{downexit}\end{equation}
so the proof consists in showing that, almost surely,  
\begin{equation}\lim_{n\rightarrow +\infty}\tau_n=+\infty .\label{downexit-diverg}\end{equation}
We will do this by means of a modification of the 
Lyapunov function used in Section 2.2 of \cite{VW85}. 

\begin{theorem}\label{off-exist+unq}
Let $W$ be a standard Brownian 
motion. Then, for every $x^0\in\overline D-\{0\}$, there is a unique 
strong solution to (\ref{SDER}) with $X(0)=x^0$. 
\end{theorem}

\begin{proof}
As anticipated above, by \cite{DI93} there is one and only 
one stochastic process $X$ that satisfies 
$(\ref{SDER})$ for $t<\lim_{n\rightarrow +\infty}\tau_n$, where $
\tau_n$ is defined 
by $(\ref{downexit})$. Therefore, we only have to prove 
$(\ref{downexit-diverg})$.
Define  
\[V(x):=|(\sigma\sigma^T)(0)^{-1/2}x|^{-p}\cos(\vartheta ((\sigma
\sigma^T)(0)^{-1/2}x)+\xi ),\]
where $\vartheta (z)\in (-\pi ,\pi ]$ is the angular polar coordinate of $
z$ 
and
\[\xi :=\vartheta ((\sigma\sigma^T)(0)^{1/2}e^{*})-2\vartheta_0,\quad
\vartheta_0:=\lim_{x\in\overline D-\{0\},\,x\rightarrow 0}\vartheta 
((\sigma\sigma^T)(0)^{-1/2}x),\quad p\in (0,1),\]
(notice that $-\frac {\pi}2<\vartheta_0<\frac {\pi}2$). Then one can check that, 
if $p$ is taken sufficiently close to $1$, 
\begin{eqnarray*}
\lim_{x\in\overline D-\{0\},x\rightarrow 0}V(x)&=&+\infty\nonumber\\
\lim_{x\in\overline D-\{0\},x\rightarrow 0}AV(x)&=&-\infty ,\nonumber\\
\lim_{x\in\partial D-\{0\},x\rightarrow 0}DV(x)\gamma (x)&=&-\infty 
.\nonumber\end{eqnarray*}
Therefore there exists $\delta >0$ such that 
\[\inf_{x\in\bar {D}-0,x_1\leq\delta}V(x)>0,\]
\[\sup_{x\in\bar {D}-0,x_1\leq\delta}AV(x)<0\]
and 
\[\sup_{x\in\partial D-0,x_1\leq\delta}DV(x)\gamma (x)<0.\]
Let
\[\alpha_{\delta}=\inf\{t\geq 0:X_1(t)\leq\delta /2\}\]
 and
\[\beta_{\delta}=\inf\{t\geq\alpha_{\delta}:X_1(t)\geq
\delta \}.\]
Without loss of  generality, we can assume that 
$x^0_1>\delta$. 

By It\^o's formula, for $n^{-1}<\delta /2$,  
\[E[V(X(\tau_n\wedge\beta_{\delta}))\bold1_{\{\alpha_{\delta}<\infty 
\}}]\leq E[V(X(\alpha_{\delta})\bold1_{\{\alpha_{\delta}<\infty \}}
].\]
Consequently, 
\[V(\frac 1n)P\{\tau_n<\beta_{\delta}|\alpha_{\delta}<\infty \}+V
(\delta )P\{\beta_{\delta}<\tau_n|\alpha_{\delta}<\infty \}\leq V
(\frac {\delta}2).\]
Consequently, if $X_1$ hits $\delta /2$, then with probability one, it 
hits $\delta$ before it hits $0$.  In particular, with 
probability one, $X_1$ never hits $0$. 
\end{proof}

\begin{remark}\label{init-distr}
Theorem \ref{off-exist+unq} implies 
existence and uniqueness of a strong solution to $(\ref{SDER})$ 
for every initial condition such that $\P (X(0)\in\overline D-\{0
\})=1$
which in turn implies  
existence and uniqueness in distribution of a weak solution to 
$(\ref{SDER})$ for every initial distribution $\mu$ such that 
$\mu (\overline D-\{0\})=1$. 
\end{remark}

\setcounter{equation}{0}

\section{Weak existence and uniqueness starting at 
$x^0=0$}\label{0}

In this section we prove weak existence and uniqueness 
for the solution of $(\ref{SDER})$ starting at the origin.

In order to prove existence (Theorem \ref{tip-exist}), 
we start with a sequence of 
solutions to $(\ref{SDER})$ starting at $x^n\in\overline D-\{0\}$, where 
$\{x^n\}$ converges to the origin. For every $n$, we consider 
a random time change of the solution, 
the same time change that is 
used in \cite{Kur90} to construct a solution to a 
patchwork martingale problem from a solution to the 
corresponding constrained martingale problem.  
The time changed processes and the time changes are relatively 
compact, and any limit point satisfies the time changed 
version of $(\ref{SDER})$ with $X(0)=0$. 
The key point of the proof is to show 
that the limit time change is invertible. 
The process obtained 
is a weak solution to 
$(\ref{SDER})$ defined for all times. 

Weak uniqueness of the solution of $(\ref{SDER})$ 
starting at the origin (Theorem \ref{tip-unq} below) 
is the main result of this paper. 
Our proof takes inspiration from the one used 
in \cite{TW93} for reflecting Brownian motion in 
the nonnegative orthant. The argument of that paper, in 
the case when the origin is not reached, 
can essentially be reformulated as follows:  
First, it is shown that, for any solution of the SDER 
starting at the origin, the 
exit time from $B_{\delta}(0)$, $\delta >0$, is finite and tends to zero as 
$\delta\rightarrow 0$, almost surely,  
and that any two solutions of the SDER, 
starting at the origin, that have the same exit 
distributions from $B_{\delta}(0)$, for all $\delta >0$ sufficiently small, 
have the same distribution; next it is proved that, for 
any $\xi\in\partial B_1(0)$, $\xi$ in the nonnegative orthant, letting $
X^{\delta\xi}$ be 
the solution of the SDER starting at $\delta\xi$ and $\tau_{2\delta}$ be its 
exit time from $B_{2\delta}(0)$, $\P\left(X^{\delta\xi}(\tau_{2\delta}
)/(2\delta )\in\cdot\right)$ is 
independent of $\delta$ and hence defines the transition kernel of a 
Markov chain on $\partial B_1(0)$. This Markov chain is shown to be ergodic 
and that in turn ensures that, for any initial distribution 
$\mu^n$ on $\overline D\cap\left(\partial B_{\delta /2^n}(0)\right
)$, the exit 
distribution of $X^{\mu^n}$ from $B_{\delta}(0)$ converges, as $n$
goes to infinity, to a uniquely determined distribution. 
Consequently, any two solutions of the SDER 
starting at the origin have the same exit 
distributions from $B_{\delta}(0)$. 

The first part of our argument is the same 
as in \cite{TW93}, except that we find it more convenient to 
use the exit distribution from $\{x:\,x_1<\delta \}$ rather than 
from $B_{\delta}(0)$. We prove that, for any solution of 
$(\ref{SDER})$ starting at the origin, the 
exit time from $\{x:\,x_1<\delta \}$, $\delta >0$, is finite and tends to zero as 
$\delta\rightarrow 0$, almost surely 
(Lemma \ref{leave}), and that any two solutions of 
$(\ref{SDER})$ starting at the origin that have the same exit 
distributions from $\{x:\,x_1<\delta \}$, for all $\delta >0$ sufficiently small, 
have the same distribution (Lemma \ref{unqexit-unq}). 

The second part of our argument consists in  
showing that, for $\{\delta_n\}$, a sequence of positive numbers 
decreasing to zero, 
any two 
solutions $X,\tilde {X}$ of (\ref{SDER}) starting at the origin satisfy 
\[{\cal L}(X_2(\tau^X_{\delta_n}))={\cal L}(\tilde {X}_2(\tau^{\tilde {
X}}_{\delta_n})),\]
where $\tau^X_{\delta_n}$ and $\tau^{\tilde {X}}_{\delta_n}$ are the corresponding exit times from 
$[0,\delta_n)$. This fact cannot be proved by the arguments used in 
\cite{TW93}. Instead, it is achieved by the rescaling 
result of Section \ref{tipscale}, together with a coupling 
argument based on \cite{LR86}.

\subsection{Existence}\label{tipexist}

\begin{theorem}\label{tip-exist}
There exists a weak solution to $(\ref{SDER})$ starting 
at $x^0=0$.
\end{theorem}

\begin{proof}
Consider a sequence $\{x^n\}\subseteq\overline D-\{0\}$ that converges to the 
origin. Let $(X^n,\Lambda^n)$ be the solution of $(\ref{SDER})$ 
starting at $x^n$. Define 
\[H^n_0(t):=\inf\{s\geq 0:\,s+\Lambda^n(s)>t\},\quad\]
and set
\[Y^n(t):=X^n(H^n_0(t)),\quad M^n(t):=W(H^n_0(t)),\quad H^n_1(t):
=\Lambda^n(H^n_0(t)),\quad\eta^n(t):=\gamma^n(H^n_0(t)).\]
Then $H^n_0,H^n_1$ are nonnegative and nondecreasing, 
\[H^n_0(t)+H^n_1(t)=t,\quad t\geq 0,\]
and 
\begin{eqnarray}
Y^n(t)=x^n+\int_0^t\sigma (Y^n(s))dM^n(s)&+&\int_0^tb(Y^n(s))dH_0^
n(s)+\int_0^t\eta^n(s)dH_1^n(s),\label{seqeq}\\
Y^n(t)\in\overline D,\quad\eta^n(t)\in\Gamma_1(Y^n(t)),\,\,\,dH_1^
n-a.e.,\quad&&\int_0^t\bold1_{\partial D}(Y^n(s))dH_1^n(s)=H_1^n(
t),\non\end{eqnarray}
where $M^n$ is a continuous, square integrable martingale with 
\[[M^n](t)=H_0^n(t)I,\quad t\geq 0.\]
With reference to Theorem 5.4 in \cite{KP91}, let
\[U^n(t)=\int_0^t\eta^n(s)dH_1^n(s).\]
Since $U_n$, $H_0^n$, and $H_1^n$ are all Lipschitz with Lipschitz 
constants bounded by $1$,
 $\left\{\left(Y^n,U^n,M^n,H^n_0,H^n_1\right)\right\}$ is relatively 
compact in distribution in the appropriate space of 
continuous function.
Taking a convergent subsequence 
with limit $\left(Y,U,M,H_0,H_1\right)$, $Y$ satisfies
\[Y(t)=Y(0)+\int_0^t\sigma (Y(s))dM(s)+\int_0^tb(Y(s))dH_0(s)+U(t
),\]
where $M(t)=W(H_0(t))$ for a standard Brownian motion 
$W$.  Since $|U^n(r)-U^n(t)|\leq |H^n_1(r)-H_1^n(r)|$, the same inequality holds for 
$U$ and $H_1$ and hence 
\[U(t)=\int_0^t\eta (s)dH_1(s).\]
It remains only to characterize $\eta$.

Invoking the Skorohod representation theorem, we assume 
that $(Y^n,U^n,M^n,H_0^n,H_1^n)\rightarrow (Y,U,M,H_0,H_1)$ uniformly over compact 
time intervals, almost surely.
Then the argument of Theorem 3.1 of \cite{Cos92} yields 
that  
\begin{equation}\eta (t)\in\Gamma_1(Y(t)),\,\,\,dH_1-a.e.,\qquad\int_
0^t1_{\partial D}(Y(s))dH_1(s)=H_1(t).\label{H1b}\end{equation}

Finally, let us show that $H_0$ is invertible and $H_0^{-1}$ is 
defined on all $[0,\infty )$. 
Suppose, by contradiction, that $H_0$ is constant 
on some time interval $[t_1,t_2]$, $0\leq t_1<t_2$. Then, since
\begin{equation}H_0(t)+H_1(t)=t.\label{H0+H1}\end{equation}
\begin{equation}H_1(s)-H_1(t)=s-t,\quad\text{{\rm f}{\rm o}{\rm r}}
t_1\leq t<s\leq t_2.\label{H0-const}\end{equation}
In particular, by $(\ref{H1b})$,  $Y(t)\in\partial D$ for all $t\in 
[t_1,t_2)$.
For $x\in\partial D-\{0\}$, let $\nu (x)$ denote the unit inward normal 
vector at $x$. If, for some $t\in [t_1,t_2)$, $Y(t)\in\partial D-
\{0\}$, then for 
$s>t$ close enough to $t$ so that $Y(r)\in\partial D-\{0\}$, $r\in 
[t,s]$, 
we have, by Condition \ref{dircond}(a),
\[\inf_{r\in [t,s]}\langle\gamma (Y(r)),\nu (Y(t))\rangle >0,\]
and hence 
\[\langle Y(s)-Y(t),\nu (Y(t))\rangle =\int_t^s\langle\gamma (Y(r
)),\nu (Y(t))\rangle dH_1(r)>0,\]
which implies, for $s$ close enough to $t$, 
\[Y(s)\in D,\]
and this contradicts $(\ref{H0-const})$. 
On the other hand, if $Y(t)=0$ for all $t\in [t_1,t_2)$, then, 
\[\int_{t_1}^t\eta (r)dH_1(r)=0,\]
while Condition \ref{dircond}(b) gives 
\[\langle\int_{t_1}^t\eta (r)dH_1(r),e^{*}\rangle\geq\inf_{\Gamma_
1(0)}\langle\gamma ,e^{*}\rangle (t-t_1)>0.\]
Therefore $H_0$ is strictly increasing. 

In order to see that 
$H_0$ diverges as $t$ goes to infinity, we can use the 
argument of Lemma 1.9 of \cite{Kur90}, provided there is 
a ${\cal C}^2$ function $\varphi$ such that 

\begin{equation}\inf_{x\in\partial D}\inf_{\gamma\in\Gamma_1(x)}D
\varphi (x)\gamma >0.\label{lemma1.9}\end{equation}
Let $e^{*}$ be the vector in Condition \ref{dircond}(b), and let 
$r^{*}>0$ be such that  
\[\inf_{x\in\partial D,\,\langle e^{*},x\rangle\leq 2r^{*}}\inf_{
\gamma\in\Gamma_1(x)}\langle e^{*},\gamma\rangle >0.\]
By Condition \ref{domcond}(b), 
\[\{x\in D:\langle e^{*},x\rangle\geq r^{*}\}=\{x\in D:\Psi_3(x)>
0,\,\langle e^{*},x\rangle\geq r^{*}\},\]
for some ${\cal C}^1$ function $\Psi_3$ such that $\inf_{x:\,\Psi_3
(x)=0}|D\Psi_3(x)|>0$. 
Then (see, e.g., \cite{CIL92}, Lemma 7.6) there exists a ${\cal C}^
2$ 
function $\varphi_3$ such that 
\[\inf_{x\in\partial D,\,\langle e^{*},x\rangle\geq r^{*}}D\varphi_
3(x)\gamma (x)>0.\]
Of course we can always assume 
\[\inf_{x\in\partial D,\,\langle e^{*},x\rangle\geq r^{*}}\varphi_
3(x)\geq 2r^{*}.\]
Therefore the function 
\[\varphi (x):=\langle e^{*},x\rangle\chi (\frac {\langle e^{*},x
\rangle -r^{*}}{r^{*}})+[1-\chi (\frac {\langle e^{*},x\rangle -r^{
*}}{r^{*}})]\varphi_3(x),\]
where $\chi :\R\rightarrow [0,1]$ is a smooth, nonincreasing function such that 
$\chi (r)=1$ for $ $$r\leq 0$, $\chi (r)=0$ for $ $$r\geq 1$, 
satisfies $(\ref{lemma1.9})$. 

We conclude our proof by setting 
\[X(t):=Y(H_0^{-1}(t)),\quad\Lambda (t):=H_1(H_0^{-1}(t)),\quad\gamma 
(t):=\eta (H_0^{-1}(t)).\]
It can be easily checked that  $X$, $\Lambda$ and $\gamma$ verify 
$(\ref{SDER})$ with $x^0=0$.
\end{proof}

\subsection{Scaling near the tip}\label{tipscale}

The following scaling result is central to our argument.

Recall Condition \ref{domcond}. Let
 $\{\delta_n\}$ be the sequence of positive numbers defined by 
\begin{equation} q_n:=(\psi_2-\psi_1)(\delta_n),\quad\delta_{n+1}:=\delta_n-q_n,
\qquad n\geq 0.\label{scale}\end{equation}
$\{\delta_n\}$ is decreasing and converges to zero. In addition, by 
Condition \ref{domcond}(c), 
\begin{equation}\lim_{n\rightarrow\infty}\frac {q_n}{\delta_n}=0.\label{scale-rate}\end{equation}

\begin{lemma}\label{scaling} 
Let $X^n$ be a solution to $(\ref{SDER})$ starting at 
$x^n\in\overline D-\{0\}$, where $\left\{\bar x^n\right\}:=\left\{
q_n^{-1}(x_1^n-\delta_{n+1},x_2^n)\right\}$ converges to 
a point $\bar {x}^0$.  

Then the sequence of processes 
\begin{equation}\left\{\bar X^n\right\}:=\left\{q_n^{-1}\left(X^n_1(q_n^2\cdot 
)-\delta_{n+1},X^n_2(q_n^2\cdot )\right)\right\}\label{scaled}\end{equation}
converges in distribution to a reflecting Brownian 
motion in $(-\infty ,\infty )\times [L,L+1]$ with directions of reflection 
$\gamma^1(0)$ on $(-\infty ,\infty )\times \{L\}$ and $\gamma^2(0
)$ on $(-\infty ,\infty )\times \{L+1\}$, 
respectively, covariance matrix $(\sigma\sigma^T)(0)$ and initial condition 
$\bar {x}^0$. 
\end{lemma}

\begin{proof}
$\bar {X}^n$ is a solution of the rescaled SDER 
\begin{eqnarray}
\bar {X}^n(t)&=&\bar {x}^n+q_n\int_0^tb((q_n\bar {X}^n_1(s)+\delta_{
n+1},q_n\bar {X}^n_2(s)))d(s)\label{scaledSDER}\\
&&\qquad\qquad +\int_0^t\sigma ((q_n\bar {X}^n_1(s)+\delta_{n+1},
q_n\bar {X}^n_2(s)))dW(s)\nonumber\\
&&\qquad\qquad +\int_0^t\bar{\gamma}^n((s)d\bar{\Lambda}^n(s),\quad 
t\geq 0,\nonumber\end{eqnarray}
and 
\[\bar {X}^n(t)\in\bar {D}^n,\quad\bar{\gamma}^n(t)\in\Gamma_1(\bar {
X}^n(t)),\,\,\,d\bar{\Lambda}^n{}^n-a.e.,\quad\bar{\Lambda}^n(t)=
\int_0^t\bold1_{\partial\bar {D}_n}(\bar {X}^n(s))d\bar{\Lambda}^
n(s),\quad t\geq 0,\]
where 
\[\bar {D}^n:=\{x:\,(q_nx_1+\delta_{n+1},q_nx_2)\in D\},\qquad\partial
\bar {D}^n:=\{x:\,(q_nx_1+\delta_{n+1},q_nx_2)\in\partial D\}.\]

Observe that 
\[\bar {D}_n\rightarrow\Delta :=(-\infty ,\infty )\times [L,L+1]\]
in the sense of 
Hausdorff convergence of sets, or, more precisely, the 
boundaries converge uniformly on compact subsets of 
$(-\infty ,\infty )$, that is, 
\[\lim_{n\rightarrow\infty}\frac {\psi_1(q_nx_1+\delta_{n+1})}{q_
n}=\lim_{n\rightarrow\infty}\frac {\psi_1(q_n(x_1-1)+\delta_n)}{q_
n}=\lim_{n\rightarrow\infty}\frac {\psi_1(\delta_n)}{q_n}=L,\]
(notice that eventually $0<q_n(x_1-1)+\delta_n<\delta_0$), and analogously 
\[\lim_{n\rightarrow\infty}\frac {\psi_2(q_nx_1+\delta_{n+1})}{q_
n}=\lim_{n\rightarrow\infty}\frac {\psi_2(\delta_n)}{q_n}=L+1.\]

Note that the second term on the right of 
(\ref{scaledSDER}) converges to zero. 
By applying the same time-change argument as in Theorem 
\ref{tip-exist}, we see that $\left\{\bar X^n\right\}$ is relatively compact 
and 
\[(q_n\bar {X}^n_1(s)+\delta_{n+1},q_n\bar {X}^n_2(s))\rightarrow 
0.\]
Consequently, $\left\{\bar X^n\right\}$ converges in distribution to $
\bar {X}$ satisfying
\begin{equation}\bar {X}(t)=\bar {x}^0+\sigma (0)W(t)+\gamma_1(0)\Lambda_L(t)+\gamma_
2(0)\Lambda_{L+1}(t),\label{reflBM}\end{equation}
where $\bar {X}(t)\in\Delta$, $\Lambda_L$ is nondecreasing and increases only 
when $\bar {X}_2=L$, and $\Lambda_{L+1}$ is nondecreasing and increases 
only when $\bar {X}_2=L+1$, that is, 
 $\bar {X}$ is a reflecting Brownian motion in 
$(-\infty ,\infty )\times [L,L+1]$ with directions of reflection 
$\gamma^1(0)$ on $(-\infty ,\infty )\times \{L\}$ and $\gamma^2(0
)$ on $(-\infty ,\infty )\times \{L+1\}$, 
respectively, covariance matrix $(\sigma\sigma^T)(0)$ and initial condition 
$\bar {x}^0$. Then the thesis 
follows from the fact that the distribution of 
$\bar {X}$ is uniquely determined. 
\end{proof}

\subsection{Uniqueness}\label{tipunq}

Of course, we can suppose, without 
loss of generality, that 
\begin{equation}\sup_{0\leq x_1\leq\delta_0}|\psi_1'(x_1)|,\,\sup_{
0\leq x_1\leq\delta_0}|\psi_2'(x_1)|<\frac 12,\label{D-RH3}\end{equation}
and that for $x_1\leq\delta_0$, $(\sigma\sigma^T)(x)$ is strictly positive definite and 
\begin{equation}\sup_{x\in\overline D,\,x_1\leq\delta_0}|\sigma (
x)^{-1}|<2|\sigma (0)^{-1}|,\quad\sup_{x\in\overline D,\,x_1\leq\delta_
0}|\sigma (x)-\sigma (0)|<\frac 12|\sigma (0)^{-1}|^{-1}.\label{D-RH4}\end{equation}

For every solution $X$ of (\ref{SDER}), let 
\begin{equation}\tau^X_{\delta}=\inf\{t\geq 0:\,X_1(t)\geq\delta \},\quad\delta 
>0.\label{exit-time}\end{equation}

\begin{lemma}\label{leave}
For $\delta$ sufficiently small, for any solution, $X$, of  
 $(\ref{SDER})$ starting at the origin 
\[\E[\tau^X_{\delta}]\leq C\delta^2.\]
\end{lemma}

\begin{proof}
Here $X$ is fixed, so we will omit the superscript $X$. 
Let $e^{*}$ be the vector in Condition \ref{dircond}(b) and 
\[f(x)=\frac 12\langle e^{*},x\rangle^2.\]
Then 
\begin{eqnarray*}\E[f(X(t\wedge\tau_{\delta}))]&=&\E\bigg[\int_0^{t\wedge\tau_{
\delta}}\big(\langle e^{*},X(s)\rangle\langle e^{*},b(X(s))\rangle 
+\frac 12(e^{*})^T(\sigma\sigma^T)(X(s))e^{*}\big)ds\\
&&\qquad\qquad\qquad\quad +\int_0^{t\wedge\tau_{\delta}}\langle e^{
*},X(s)\rangle\langle e^{*},\gamma (s)\rangle d\Lambda (s)\bigg].\end{eqnarray*}
Observe that 
\[\lim_{x\in\overline D-\{0\},\,x_1\rightarrow 0^{+}}\frac {|x|}{
x_1}=1,\qquad\lim_{x\in\overline D-\{0\},\,x_1\rightarrow 0^{+}}\frac {
\langle e^{*},x\rangle}{e^{*}_1x_1}=1,\]
and hence, for $\delta$ sufficiently small, 
\[|X(t\wedge\tau_{\delta})|^2\leq 4\delta^2,\]
and 
\[\int_0^{t\wedge\tau_{\delta}}\langle e^{*},X(s)\rangle\langle e^{
*},\gamma (s)\rangle d\Lambda_i(s)\geq\frac 12\int_0^{t\wedge\tau_{
\delta}}e^{*}_1X_1(s)\langle e^{*},\gamma (s)\rangle d\Lambda (s)
.\]
Therefore, for $\delta$ sufficiently small, 
\[2\delta^2\geq\E[f(X(t\wedge\tau_{\delta}))]\geq\frac 14(e^{*})^
T(\sigma\sigma^T)(0)e^{*}\E[t\wedge\tau_{\delta}],\]
which yields the assertion by taking the limit as $t$ goes to 
infinity.
\end{proof}

\begin{remark}\label{leavex0}
By looking at the proof of Lemma \ref{leave}, we see that 
we have proved, more generally, that, for every $x^0\in\overline 
D$ 
with $x^0_1<\delta$, for every solution $X$ of $(\ref{SDER})$ starting at $
x^0$, 
\[\E[\tau^X_{\delta}]\leq C\left(4\delta^2-\langle e^{*},x^0\rangle^
2\right).\]
\end{remark}

\begin{lemma}\label{unqexit-unq}
Suppose any two weak solutions, $X,\tilde {X}$, of $(\ref{SDER})$ starting 
at the origin satisfy  
\begin{equation}{\cal L}(X(\tau^X_{\delta}))={\cal L}(\tilde {X}(\tau^{\tilde {X}}_{
\delta})),\label{eqdist}\end{equation}
for all $\delta$ sufficiently small (recall that, by Lemma \ref{leave}, $
\tau^X_{\delta}$ 
and $\tau^{\tilde {X}}_{\delta}$ are almost surely finite). Then the solution of 
$(\ref{SDER})$ starting at the origin is unique in 
distribution. 
\end{lemma}

\begin{proof}
Since starting away from $0$, we have strong and weak uniqueness, 
if (\ref{eqdist}) holds, we have ${\cal L}(X(\tau_{\delta}^X+\cdot 
))={\cal L}(\tilde {X}(\tau^{\tilde {X}}_{\delta}+\cdot )$.  
Since  $\tau^X_{\delta}$ and $\tau^{\tilde {X}}_{\delta}$ converge to zero as $
\delta\rightarrow 0$ and 
$X(\tau_{\delta}^X+\cdot )\rightarrow X$ and $\tilde {X}(\tau^{\tilde {
X}}_{\delta}+\cdot )\rightarrow\tilde {X}$, uniformly over compact 
time intervals.
Consequently, we must have ${\cal L}(X)={\cal L}(\tilde {X})$. 
\end{proof}

\begin{theorem}\label{tip-unq}
The solution of 
$(\ref{SDER})$ starting at $x^0=0$ is unique in 
distribution.
\end{theorem}

\begin{proof}
In what follows, $X$ and $\tilde {X}$ will be weak solutions of $
(\ref{SDER})$ 
starting at the origin.  Let  
$\{\delta_n\}=\{\delta_n\}_{n\geq 0}$ be given by (\ref{scale}). 
We want to show that 
\begin{equation}\Vert{\cal L}(X_2(\tau^X_{\delta_n}))-{\cal L}(\tilde {
X}_2(\tau^{\tilde {X}}_{\delta_n}))\Vert_{TV}\leq (1-p_0\eta^2_0)
\Vert{\cal L}(X_2(\tau^X_{\delta_{n+2}}))-{\cal L}(\tilde {X}_2(\tau^{
\tilde {X}}_{\delta_{n+2}}))\Vert_{TV},\label{tvsol}\end{equation}
where $p_0,\eta_0\in (0,1)$ come from Lemma \ref{coupling} and Lemma 
\ref{positive-hitting}, respectively, so that 
by iterating (\ref{tvsol}), we obtain 
\[\Vert {\cal L}(X_2(\tau^X_{\delta_n}))-{\cal L}(\tilde {X}_2(\tau^{
\tilde {X}}_{\delta_n}))\Vert_{TV}=0,\]
for all $n$.  The theorem then follows from Lemma 
\ref{unqexit-unq}.

To prove this, we construct below two solutions of (\ref{SDER}), 
$\chi$, starting at $(\delta_{n+2},\chi_2(0))$, and $\tilde{\chi}$, starting at $
(\delta_{n+2},\tilde{\chi}_2(0))$, that are 
coupled in such a way that, letting 
\begin{equation}\tau_2:=\inf\{t\geq 0:\,\chi_1(t)\geq\delta_n\},\quad
\tilde{\tau}_2:=\inf\{t\geq 0:\,\tilde{\chi}_1(t)\geq\delta_n\},\label{chi-ex}\end{equation}
it holds
\[{\Bbb P}\big(\chi (\tau_2)=\tilde{\chi }(\tilde{\tau}_2)\big)\geq 
p_0\eta_0^2,\]
no matter what the distributions of $\chi_2(0)$ and $\tilde{\chi}_
2(0)$. 
This implies that 
\[\Vert {\cal L}(\chi (\tau_2))-{\cal L}(\tilde{\chi }(\tilde{\tau}_
2))\Vert_{TV}\leq (1-p_0\eta_0^2).\]
Consequently, with the notation of Lemma \ref{tv}, 
denoting by $P$ the transition function from 
$[\psi_1(\delta_{n+2}),\psi_2(\delta_{n+2})]$ to $[\psi_1(\delta_
n),\psi_2(\delta_n)]$ defined by 
\[P(x_2,\cdot ):=\P\big(\chi_2(\tau_2)\in\cdot\big|\chi_2(0)=x_2\big
)=\P\big(X_2(\tau^X_{\delta_n})\in\cdot\big|X_2(\tau^X_{\delta_{n
+2}})=x_2\big),\]
we have 
\[\Vert P\nu -P\tilde{\nu}\Vert_{TV}\leq (1-p_0\eta_0^2),\]
for any two probability distributions $\nu$ and $\tilde{\nu}$ on 
$[\psi_1(\delta_{n+2}),\psi_2(\delta_{n+2})]$. Therefore (\ref{tvsol}) follows from 
Lemma \ref{tv}.

We conclude the proof with the construction of the 
coupled solutions $\chi$ and $\tilde{\chi}$. We start $\chi$ and $
\tilde{\chi}$ as two 
independent solutions of $(\ref{SDER})$, with initial 
condition $(\delta_{n+2},\chi_2(0))$ and $(\delta_{n+2},\tilde{\chi}_
2(0))$, respectively, 
and we run them until the times 
\[\tau_1:=\inf\{t\geq 0:\,\chi_1(t)\geq\delta_{n+1}\},\quad\tilde{
\tau}_1:=\inf\{t\geq 0:\,\tilde{\chi}_1(t)\geq\delta_{n+1}\}.\]
We then consider a solution, $(Z,\tilde {Z})$, 
with initial distribution $(\chi (\tau_1),\tilde{\chi }(\tau_1))$, 
of the coupled SDE $(\ref{coupleq})$ with $\beta =b$, $\varsigma 
=\sigma$ and $B$ independent of 
$(\chi (\tau_1),\tilde{\chi }(\tau_1))$, until the times 
\[\Theta :=\inf\{t\geq 0:\,Z(t)\notin Q\},\quad\tilde{\Theta }:=\inf
\{t\geq 0:\,\tilde {Z}(t)\notin Q\},\]
where $Q$ is the rectangle 
\[Q:=(\delta_{n+1}-\frac 14q_{n+1},\delta_{n+1}+\frac 14q_{n+1})\times 
I^{\epsilon_0+1/2}_{n+1},\]
and $I^{\epsilon_0+1/2}_{n+1}$ is the interval in Lemma \ref{positive-hitting} for 
a value of $\epsilon_0$ to be chosen later. We set 
\[\chi (\tau_1+t):=Z(t),\quad\tilde{\chi }(\tau_1+t):=\tilde {Z}(
t),\qquad\mbox{\rm for }t\leq\Theta\wedge\tilde{\Theta }.\]
For $\epsilon_0\leq 1/4$ , by $(\ref{D-RH3})$, 
\[\overline Q\subseteq D,\]
therefore $\chi$ and $\tilde{\chi}$ are solutions of $(\ref{SDER}
)$ up to 
$\tau_1+\Theta\wedge\tilde{\Theta}$  and $\tilde{\tau}_1+\Theta\wedge
\tilde{\Theta}$ respectively. 
Moreover, by $(\ref{D-RH4})$, the assumptions of Lemma 
\ref{coupling} are satisfied. By choosing $\epsilon_0\leq C_0/4$, 
where $C_0$ is the constant in Lemma \ref{coupling}, we 
have that, for $x_2,\,\tilde {x}_2\in I^{\epsilon_0}_{n+1}$, it holds 
\[|x_2-\tilde {x}_2|\leq C_0\frac {q_{n+1}}4,\qquad\overline {B_{
q_{n+1}/4}(\delta_{n+1},x_2)},\,\overline {B_{q_{n+1}/4}(\delta_{
n+1},\tilde {x}_2)}\subseteq\overline Q.\]
Therefore Lemma \ref{coupling} yields that 
\[\P\big(Z(\Theta\wedge\tilde{\Theta })=\tilde {Z}(\Theta\wedge\tilde{
\Theta })\big|Z_2(0)\in I^{\epsilon_0}_{n+1},\tilde {Z}_2(0)\in I^{
\epsilon_0}_{n+1}\big)\geq p_0.\]
Combining this with Lemma \ref{positive-hitting}, we 
get 
\begin{eqnarray}
&&\P\big(\chi (\tau_1+\Theta\wedge\tilde{\Theta })=\tilde{\chi }(
\tilde{\tau}_1+\Theta\wedge\tilde{\Theta })\big)\non\\
&\geq&\P\big(\chi (\tau_1+\Theta\wedge\tilde{\Theta })=\tilde{\chi }
(\tilde{\tau}_1+\Theta\wedge\tilde{\Theta })\big|\chi_2(\tau_1)\in 
I^{\epsilon_0}_{n+1},\tilde{\chi}_2(\tilde{\tau}_1)\in I^{\epsilon_
0}_{n+1}\big)\,\P\big(\chi_2(\tau_1)\in I^{\epsilon_0}_{n+1},\tilde{
\chi}_2(\tilde{\tau}_1)\in I^{\epsilon_0}_{n+1}\big)\non\\
&\geq&\P\big(Z(\Theta\wedge\tilde{\Theta })=\tilde {Z}(\Theta\wedge
\tilde{\Theta })\big|Z_2(0)\in I^{\epsilon_0}_{n+1},\tilde {Z}_2(
0)\in I^{\epsilon_0}_{n+1}\big)\,\eta_{\epsilon_0}^2\non\\
&\geq&p_0\eta_{\epsilon_0}^2.\non\end{eqnarray}
Finally, we define $\chi (\tau_1+\Theta\wedge\tilde{\Theta }+\cdot 
)$ as a solution of 
$(\ref{SDER})$ starting at $\chi (\tau_1+\Theta\wedge\tilde{\Theta }
)$ and 
$\tilde{\chi }(\tilde{\tau}_1+\Theta\wedge\tilde{\Theta }+\cdot )
=\chi (\tau_1+\Theta\wedge\tilde{\Theta }+\cdot )$, if 
$\chi (\tau_1+\Theta\wedge\tilde{\Theta })=\tilde{\chi }(\tilde{\tau}_
1+\Theta\wedge\tilde{\Theta })$, and as a  
solution of $(\ref{SDER})$ starting at $\tilde{\chi }(\tilde{\tau}_
1+\Theta\wedge\tilde{\Theta })$ 
otherwise. Since, with $Q$ as above, 
\[\tau_1+\Theta\wedge\tilde{\Theta }<\tau_2,\quad\tilde{\tau}_1+\Theta
\wedge\tilde{\Theta }<\tilde{\tau}_2,\]
$\chi$ and $\tilde{\chi}$ have the desired property.
\end{proof}

\subsection{The Feller property}\label{Feller}

We conclude with the observation that the family of 
distributions $\big\{P^x\big\}_{x\in\overline D}$, where $P^x$ is the distribution of 
the unique weak solution of (\ref{SDER}) starting at $x$, 
enjoys the Feller property.

\begin{proposition}\label{initcont} 
Let $X^x$ be the unique weak solution of (\ref{SDER}) starting 
at $x$. Then the mapping 
$x\in\overline D\rightarrow X^x$ is continuous in distribution. 
\end{proposition}

\begin{proof}
The proof is exactly the same as that of Theorem 
\ref{tip-exist}. In fact, once it is known that the weak 
solution of (\ref{SDER}) starting at the origin is unique, 
the proof of Lemma \ref{tip-exist} amounts to showing 
that $X^x$ is continuous in distribution at the origin. 
\end{proof}

\setcounter{equation}{0}

\section{Technical lemmas}\label{technic}

\begin{lemma}
Let $\bar {X}^{(0,x_2)}$ satisfy (\ref{reflBM}) with $\bar {X}^{(
0,x_2)}(0)=(0,x_2)$, 
$x_2\in [L,L+1]$, and let
\[\bar{\tau}_1:=\inf\{t\geq 0:\,\bar {X}_1^{(0,x_2)}(t)\geq 1\}.\]
Then
\begin{itemize}
\item[(i)] $\bar{\tau}_1$ is a.s. finite.
\item[(ii)] For every $\epsilon >0$, 
\begin{equation}\inf_{L\leq x_2\leq L+1}\P\big(\,\,\bar {X}^{(0,x_
2)}_2(\bar{\tau}_1)\in (L+\frac 12-\epsilon ,L+\frac 12+\epsilon 
)\big)>0.\label{BMpositivehit}\end{equation}
\end{itemize}
\end{lemma}

\begin{proof}

\noindent{\bf (i)} To simplify notation, whenever possible 
without loss of clarity, we will omit the 
superscript on $\bar {X}$. 

Let $e^{*}$ be the vector in Condition \ref{dircond}(b). Then, 
for $\bar {X}(0)=x^0\in\R\times [L,L+1]$, for all 
$N>0$, 
\[\inf\{t\geq 0:\,\langle\bar {X}(t),e^{*}\rangle\geq N+\langle x^
0,e^{*}\rangle \}\leq\inf\{t\geq 0:\,\langle\sigma (0)W(t),e^{*}\rangle
\geq N\}<+\infty\quad a.s..\]
On the other hand 
\[\bar {X}_1(t)=\frac 1{e_1^{*}}\big(\langle\bar {X}(t),e^{*}\rangle 
-\bar {X}_2(t)e_2^{*}\big)\geq\frac 1{e_1^{*}}\big(\langle\bar {X}
(t),e^{*}\rangle -(|L|+1)\big).\]
Therefore, for $N$ large enough, 
\[\bar{\tau}_1\leq\inf\{t\geq 0:\,\langle\bar {X}(t),e^{*}\rangle
\geq N+\langle x^0,e^{*}\rangle \}.\]

\noindent{\bf (ii)} We can suppose, without loss of generality, 
$\epsilon <1/2$.  Let  
$h:\R^2\rightarrow [0,1]$ be a smooth function such that 
\[h(1,L+\frac 12)=1,\qquad h(x)=0\mbox{\rm \ for }x\notin B_{\epsilon}
((1,L+\frac 12)).\]
We can estimate the probability in the left hand side of 
$(\ref{BMpositivehit})$ by 
\[
\P\big(\bar {X}_2(\bar{\tau}_1)\in (L+\frac 12-\epsilon ,L+\frac 
12+\epsilon )\big)\geq\E\big[h(\bar {X}(\bar{\tau}_1))\big].\]
Set  
\[u(x):=\E\big[h(\bar {X}^x(\bar{\tau}_1^x))\big],\quad x\in (-\infty 
,1]\times [L,L+1].\]
$u$ is continuous on $(-\infty ,1]\times [L,L+1]$, because $\bar {
X}$ is a Feller 
process and the functional: $\inf\{t\geq 0:\,\,x(t)\geq 1\}$, 
$x(\cdot )\in {\cal C}_{\R\times [L,L+1]}[0,\infty )$, is 
almost surely continuous under the law of $\bar {X}^x$, for every 
$x\in (-\infty ,1]\times [L,L+1]$. For any bounded, smooth domain $
Q$ 
such that $\overline Q\subseteq (-\infty ,1)\times (L,L+1)$, $u$ is the classical solution 
of the Dirichlet problem with itself as boundary datum.  
Therefore $u\in {\cal C}^2((-\infty ,1)\times (L,L+1))$ and 
\[\mbox{\rm tr}((\sigma\sigma^T)(0)D^2u(x))=0,\quad\forall x\in (
-\infty ,1)\times (L,L+1).\]
For $0<\eta <1/2$, let
\[Q_{\eta}:=(-1+\eta ,1-\eta )\times (L+\eta ,L+1-\eta ).\]
By the Harnack inequality, 
\begin{equation}\inf_{x\in Q_{\eta}}u(x)\geq c_{\eta}\sup_{x\in Q_{
\eta}}u(x),\label{Harnack}\end{equation}
for some $c_{\eta}>0$. For $\eta$ small enough the right hand side of $
(\ref{Harnack})$ 
is strictly positive and hence 
\[u(0,x_2)>0,\qquad\forall x_2\in (L,L+1).\]
Now let 
\[\bar{\tau}_{Q_{\eta}}^{(0,L+1)}:=\inf\{t\geq 0:\,\bar {X}^{(0,L
+1)}(t)\in\overline {Q_{\eta}}\},\]
and fix $\eta$ small enough that the right hand side of 
$(\ref{Harnack})$ is strictly positive and that 
\[\P\big(\bar{\tau}_{Q_{\eta}}^{(0,L+1)}<\bar{\tau}_1^{(0,L+1)})\big
)>0.\]
Then 
\begin{eqnarray*}
u(0,L+1)&=&\E\big[h(\bar {X}^{(0,L+1)}(\bar{\tau}_1^{(0,L+1)}))\big
]\\
&\geq&\E\big[\I_{\{\bar{\tau}_{Q_{\eta}}^{(0,L+1)}<\bar{\tau}_1^{
(0,L+1)})\}}h(\bar {X}^{(0,L+1)}(\bar{\tau}_1^{(0,L+1)}))\big]\\
&=&\E\big[\I_{\{\bar{\tau}_{Q_{\eta}}^{(0,L+1)}<\bar{\tau}_1^{(0,
L+1)})\}}u(\bar {X}^{(0,L+1)}(\bar{\tau}_{Q_{\eta}}^{(0,L+1)}))\big
]\\
&\geq&\inf_{x\in Q_{\eta}}u(x)\,\,\,\P\big(\bar{\tau}_{Q_{\eta}}^{
(0,L+1)}<\bar{\tau}_1^{(0,L+1)})\big)\\
&>&0.\end{eqnarray*}
Analogously 
\[u(0,L)>0,\]
and the assertion follows by the continuity of $u$ on 
$\{0\}\times [L,L+1]$.
\end{proof}

\begin{lemma}\label{positive-hitting}
Let $X^{x^0}$ be the solution of $(\ref{SDER})$
starting at $x^0\in\overline D-\{0\}$, $\tau_{\delta_n}^{x^0}:=\inf
\{t\geq 0:\,X_1^{x^0}(t)\geq\delta_n\}$, 
$n\geq 0$, and for $0<\epsilon <1$, $I^{\epsilon}_n$ be the open interval 
of length $\epsilon (\psi_2-\psi_1)(\delta_n)$ centered at $\frac {
(\psi_1+\psi_2)(\delta_n)}2$. 
Then there exists $n_{\epsilon}\geq 0$ and $\eta_{\epsilon}>0$ such that 
\begin{equation}\inf_{n\geq n_{\epsilon}}\,\,\,\inf_{x_2:\,(\delta_{
n+1},x_2)\in\overline D}\P\big(X_2^{(\delta_{n+1},x_2)}(\tau_{\delta_
n}^{(\delta_{n+1},x_2)})\in I^{\epsilon}_n\big)=\eta_{\epsilon}>0
.\label{poshit}\end{equation}
\end{lemma}

\begin{proof}
Let $\{x_2^n\}$, $\psi_1(\delta_{n+1})\leq x_2^n\leq\psi_2(\delta_{
n+1})$ be such that $q_n^{-1}x^n_2$ converges 
to $\bar {x}_2^0\in [L,L+1]$, and let $X^n$ denote the solution of 
$(\ref{SDER})$ starting at $(\delta_{n+1},x^n_2)$. Let
 $\bar {X}^n$ denote the scaled process 
$(\ref{scaled})$ with initial condition $(0,q_n^{-1}x_2^n)$, and let $ $$
\bar {X}$ 
denote the limiting reflecting Brownian motion starting at 
$(0,\bar {x}_2^0)$.  Define 
\begin{eqnarray*}
\tau^n=\tau_{\delta_n}^{(\delta_{n+1},x_2^n)}&:=&\inf\{t\geq 0:\,
X_1^n(t)\geq\delta_n\}\\
\bar{\tau}^n&:=&\inf\{t\geq 0:\,\bar {X}_1^n(t)\geq 1\}\\
\bar{\tau}_1&:=&\inf\{t\geq 0:\,\bar {X}_1(t)\geq 1\}\end{eqnarray*}
Notice that $\tau^n$ is a.s. finite by Remark \ref{leavex0} and that 
\[\tau^n=q_n^{-2}\bar{\tau}^n.\]
Since the first exit time from $(-\infty ,1)\times\R$ is a 
continuous functional on a set of paths that 
has probability one under the distribution of $\bar {X}$, by 
the continuous mapping theorem we may assume that 
$\bar {X}_2^n(\bar{\tau}^n)$ converges in distribution to $\bar {
X}_2(\bar{\tau}_1)$. Then 
\begin{eqnarray*}
\liminf_n\P\left(X_2^n(\tau^n)\in I^{\epsilon}_n\right)&\geq&\liminf_
n\P\left(X_2^n(\tau^n)\in (q_n(L+\frac 12-\frac {\epsilon}4),q_n(
L+\frac 12+\frac {\epsilon}4))\right)\\
&=&\liminf_n\P\left(\bar X_2^n(\bar\tau^n)\in (L+\frac 12-\frac {
\epsilon}4,L+\frac 12+\frac {\epsilon}4)\right)\\
&\geq&\P\left(\bar X_2(\bar\tau_1)\in (L+\frac 12-\frac {\epsilon}
4,L+\frac 12+\frac {\epsilon}4)\right)\\
&\geq&\inf_{L\leq x_2\leq L+1}\P\big(\,\,\bar {X}^{(0,x_2)}_2(\bar{
\tau}_1^{(0,x_2)})\in (L+\frac 12-\frac {\epsilon}4,L+\frac 12+\frac {
\epsilon}4)\big),\end{eqnarray*}
and the assertion follows by 
$(\ref{BMpositivehit})$ and by the arbitrariness of $\{x_2^n\}$.
\end{proof}

The following lemma, which uses the 
coupling of \cite{LR86}, may be of independent 
interest. 

\begin{lemma}\label{coupling}
Let $\beta :\R^d\rightarrow\R^d$ and $\varsigma :\R^d\rightarrow\R^
d\times\R^d$ be Lipschitz 
continuous and bounded and let $\varsigma\varsigma^T$ be uniformly positive 
definite. Suppose that 
\begin{equation}\begin{array}{c}
\sup_{x,\tilde {x}}|\varsigma (x)-\varsigma (\tilde {x})|<2(\sup_
x|\varsigma (x)^{-1}|)^{-1}.\end{array}
\label{q-const}\end{equation}
Define 
\begin{equation}K(x,\tilde {x}):=I-2\frac {\varsigma (\tilde {x})^{
-1}(x-\tilde {x})(x-\tilde {x})^T(\varsigma (\tilde {x})^{-1})^T}{
|\varsigma (\tilde {x})^{-1}(x-\tilde {x})|^2},\label{couple-coeff}\end{equation}
Let $B$ be a standard Brownian motion on a probability 
space $(\Omega ,{\cal F},\P )$ and $(Z,\tilde {Z})$ be the solution of the system of 
stochastic differential equations 
\begin{eqnarray}
&&dZ(t)=\beta (Z(t)dt+\varsigma (Z(t))dB(t),\quad Z(0)=x^0,\label{coupleq}\\
&&d\tilde {Z}(t)=\beta (\tilde {Z}(t)dt+\varsigma (\tilde {Z}(t))
K(Z(t),\tilde {Z}(t))dB(t),\quad\tilde {Z}(0)=\tilde {x}^0\neq x^
0\non\end{eqnarray}
for $t<\zeta$, where 
\begin{equation}\zeta :=\lim_{\epsilon\rightarrow 0}\zeta_{\epsilon}
,\quad\zeta_{\epsilon}:=\inf\{t\geq 0:\,|Z(t)-\tilde {Z}(t)|<\epsilon 
\},\label{couple-time}\end{equation}
and $\tilde {Z}(t)=Z(t)$ for $t\geq\zeta$, on the set $\{\zeta <\infty 
\}$. (Notice that 
$K$ is locally Lipschitz continuous on 
$\R^d\times\R^d-\{(x,x),\,x\in\R^d\}$.)

Then $\tilde {Z}$ is a diffusion process with generator 
$Df(x)\beta (x)+\frac 12\mbox{\rm tr}(\varsigma (x)\varsigma (x)^
TD^2f(x))$ and, for every $p_0$, 
$0<p_0<\frac 14$, there exists a positive constant $C_0<1$, 
depending only on $p_0$, $\beta$ and $\varsigma$, such that, setting 
\[\vartheta_{\rho}:=\inf\{t\geq 0:\,|Z(t)-x^0|>\rho \},\quad\tilde{
\vartheta}_{\rho}:=\inf\{t\geq 0:\,|\tilde {Z}(t)-\tilde {x}^0|>\rho 
\},\]
for $\rho\leq 1$, $ $
\[|x^0-\tilde {x}^0|\leq C_0\rho\quad\mbox{\rm implies}\quad\P (\zeta
\leq\vartheta_{\rho}\wedge\tilde{\vartheta}_{\rho})\geq p_0.\]
\end{lemma}

\begin{proof}
The fact that $\tilde {Z}$ has generator $Df(x)\beta (x)+\frac 12\mbox{\rm tr}
(\varsigma (x)\varsigma (x)^TD^2f(x))$ 
follows from the fact $K(x,\tilde {x}$) is an orthogonal matrix. 

As in \cite{LR86}), consider  
\[U(t):=|Z(t)-\tilde {Z}(t)|.\]
For $t<\zeta$, $U$ satisfies 
\[dU(t)=a(t)dt+\alpha (t)dW(t),\]
where 
\begin{eqnarray*}
a(t)&:=&\frac {\langle Z(t)-\tilde {Z}(t),\beta (Z(t))-\beta (\tilde {
Z}(t))\rangle}{|Z(t)-\tilde {Z}(t)|}+\,\frac {\mbox{\rm tr}\big(\big
(\varsigma (Z(t))-\varsigma (\tilde {Z}(t))\big)\big(\varsigma (Z
(t))-\varsigma (\tilde {Z}(t))\big)^T\big)}{2|Z(t)-\tilde {Z}(t)|}\\
&&\quad\quad -\,\frac {\big|(\varsigma (Z(t))-\varsigma (\tilde {
Z}(t)))^T\frac {Z(t)-\tilde {Z}(t)}{|Z(t)-\tilde {Z}(t)|}\big|^2}{
2|Z(t)-\tilde {Z}(t)|}\\
\alpha (t)&:=&\left|\big(\varsigma (Z(t))-\varsigma (\tilde Z(t))
K(Z(t),\tilde Z(t))\big)^T\frac {Z(t)-\tilde {Z}(t)}{|Z(t)-\tilde {
Z}(t)|}\right|,\end{eqnarray*}
and $W$ is a standard Brownian motion. 
Then, as in \cite{LR86}),  setting
\begin{eqnarray*}
g(u)&:=&\sup_{|x-\tilde {x}|=u}\bigg\{\big|\big(\varsigma (x)-\varsigma 
(\tilde {x})K(x,\tilde {x})\big)^T\frac {x-\tilde {x}}{|x-\tilde {
x}|}\big|\,^{-2}\\
&&\times\bigg[\frac {\langle x-\tilde {x},\beta (x)-\beta (\tilde {
x})\rangle}{|x-\tilde {x}|}+\frac {{\rm t}{\rm r}\big((\varsigma 
(x)-\varsigma (\tilde {x}))(\varsigma (x)-\varsigma (\tilde {x}))^
T\big)-\big|(\varsigma (x)-\varsigma (\tilde {x}))^T\frac {x-\tilde {
x}}{|x-\tilde {x}|}\big|^2}{2|x-\tilde {x}|}\bigg]\bigg\},\end{eqnarray*}
we have 
\[a(t)\leq\alpha (t)^2g(U(t)).\]
In addition, by $(\ref{q-const})$, the Lipschitz 
property of $\beta$ and $\varsigma$ and the boundedness of $\varsigma$,  
\begin{eqnarray}
\underline {\alpha}:=\inf_t\alpha (t)&>&0\nonumber\\
\overline {\alpha}:=\sup_t\alpha (t)&\leq&2\|\varsigma\|\label{g-bound}\\
|a(t)|&\leq&\overline aU(t),\mbox{\rm \quad a.s.}\nonumber\\
|g(u)|&\leq&C\,u\nonumber\end{eqnarray}
(See the computations at page 866 of \cite{LR86}.) In 
particular $g$ is locally square integrable on $(0,\infty )$. 
Therefore, applying It\^o's formula to the process $U$ and the function 
\[\begin{array}{c}
G(u):=\int_1^udv\,\exp\big(-2\int_1^vg(z)dz\big),\end{array}
\]
we see that, for $u:=|x^0-\tilde {x}^0|$, $\theta_{2u}:=\inf\{t\geq 
0:\,U(t)>2u\}$, 
\[\P (\zeta_{\epsilon}<\theta_{2u})\geq\frac {G(2u)-G(u)}{G(2u)-G
(\epsilon )}.\]

On the other hand, we have, for any $0<t\leq 1$,  
\begin{eqnarray*}
\P (\zeta_{\epsilon}<\vartheta_{\rho}\wedge\tilde{\vartheta}_{\rho}
)&\geq&\P (\zeta_{\epsilon}\leq t,\,\,\,t<\vartheta_{\rho}\wedge\tilde{
\vartheta}_{\rho})\\
&\geq&\P (\zeta_{\epsilon}\leq t)-\P (\,\vartheta_{\rho}\wedge\tilde{
\vartheta}_{\rho}\leq t)\\
&\geq&\P (\zeta_{\epsilon}<\theta_{2u},\,\,\,\zeta_{\epsilon}\wedge
\theta_{2u}\leq t)-\P (\vartheta_{\rho}\wedge\tilde{\vartheta}_{\rho}
\leq t)\\
&\geq&\P (\zeta_{\epsilon}<\theta_{2u})-\P (\zeta_{\epsilon}\wedge
\theta_{2u}>t)-\P (\vartheta_{\rho}\wedge\tilde{\vartheta}_{\rho}
\leq t)\\
&\geq&\frac {G(2u)-G(u)}{G(2u)-G(\epsilon )}-\P\big(\sup_{s\leq t}\big
|U(s)-u\big|\leq u\big)-\P (\vartheta_{\rho}\wedge\tilde{\vartheta}_{
\rho}\leq t).\end{eqnarray*}
By $(\ref{g-bound})$, we have 
\[G(0):=\lim_{\epsilon\rightarrow 0}G(\epsilon )>-\infty .\]
Then, noting that $\P (\vartheta_{\rho}\wedge\tilde{\vartheta}_{\rho}
<\infty )=1$, 
taking the limit as $\epsilon$ goes to zero, we obtain 
\begin{eqnarray*}
\P (\zeta\leq\vartheta_{\rho}\wedge\tilde{\vartheta}_{\rho})&\geq&\frac {
G(2u)-G(u)}{G(2u)-G(0)}-\P\bigg(\sup_{s\leq t}\big|\int_0^sa(r)dr
+\int_0^s\alpha (r)dW(r)\big|\leq u\bigg)\\
&&\qquad\qquad\qquad -\P (\vartheta_{\rho}\wedge\tilde{\vartheta}_{
\rho}\leq t).\end{eqnarray*}
Now, we can easily see (for instance 
applying It\^o's formula to the function $f(x)=|x-x^0|^2$) 
that  
\[\P (\vartheta_{\rho}\wedge\tilde{\vartheta}_{\rho}\leq t)\leq\frac {
C_1t}{\rho^2},\]
where $C_1$ depends only on $b$ and $\varsigma$. Of course we can 
suppose, without loss of generality, that $C_1\geq 1$. 
Therefore, we will take 
\begin{equation}t=\frac 1{C_1}(\frac 14-p_0)\,\rho^2.\label{t}\end{equation}
We have, for $\rho\leq 1$ and $u\leq C_0\rho$, where $C_0$ is a constant to be 
chosen later, 
\begin{eqnarray*}
&&\P\bigg(\sup_{s\leq t}\big|\int_0^sa(r)dr+\int_0^s\alpha (r)dW(
r)\big|\leq u\bigg)\\
&=&\P\bigg(\sup_{s\leq t}\big|\int_0^sa(r)dr+\int_0^s\alpha (r)dW
(r)\big|\leq u,\,\sup_{s\leq t}U(s)\leq 2u\bigg)\\
&\leq&\P\bigg(\sup_{s\leq t}\big|\int_0^t\alpha (s)dW(s)\big|\leq 
C_0(1+2\overline a)\sqrt {\frac {4C_1}{1-4p_0}}\,\sqrt {t}\bigg),\end{eqnarray*}
where in the last inequality we have used the fact that 
$|a(t)|\leq\overline aU(t)$, $(\ref{g-bound})$ and $(\ref{t})$. 
We take $C_0$ 
small enough that 
\begin{equation}C:=C_0(1+2\overline {\alpha})\sqrt {\frac {4C_1}{
1-4p_0}}\,<\underline {\alpha}.\label{C}\end{equation}
Then setting $\theta^{\alpha}_{C\sqrt {t}}\,:=\inf\{s\geq 0:\,\big
|\int_0^s\alpha (r)dW(r)\big|>C\sqrt {t}\}$, 
\begin{eqnarray*}
&&\P\bigg(\sup_{s\leq t}\big|\int_0^t\alpha (s)dW(s)\big|\leq C\sqrt {
t}\bigg)\\
&&\qquad\qquad =\P\bigg(\theta^{\alpha}_{C\sqrt {t}}\,\geq t,\,\,\big
(\int_0^t\alpha (s)dW(s)\big)^2\leq C^2t\bigg)\\
&&\qquad\qquad =\P\bigg(\theta^{\alpha}_{C\sqrt {t}}\,\geq t,\,\,
2\int_0^t\big(\int_0^s\alpha (r)dW(r)\big)\alpha (s)dW(s)+\int_0^
t\alpha (s)^2ds\leq C^2t\bigg)\\
&&\qquad\qquad\leq\P\bigg(\int_0^{t\wedge\theta^{\alpha}_{C\sqrt {
t}}\,}\big(\int_0^s\alpha (r)dW(r)\big)\alpha (s)dW(s)\leq\frac 1
2(C^2-\underline {\alpha}^2)t\bigg),\end{eqnarray*}
where the last inequality uses $(\ref{g-bound})$. Since 
$C<\underline {\alpha}$, the above chain of inequalities can be continued as 
\begin{eqnarray*}
&\leq&\P\bigg(\bigg|\int_0^{t\wedge\theta^{\alpha}_{C\sqrt {t}}\,}\big
(\int_0^s\alpha (r)dW(r)\big)\alpha (s)dW(s)\bigg|\geq\frac 12(\underline {
\alpha}^2-C^2)t\bigg)\\
&\leq&\frac 4{(\underline {\alpha}^2-C^2)^2\,t^2}\,\,\E\bigg[\int_
0^{t\wedge\theta^{\alpha}_{C\sqrt {t}}\,}\big(\int_0^s\alpha (r)d
W(r)\big)^2\alpha (s)^2ds\bigg]\\
&\leq&\frac {4C^2t}{(\underline {\alpha}^2-C^2)^2\,t^2}\E\bigg[\int_
0^t\alpha (s)^2ds\bigg]\\
&\leq&\frac {4C^2\bar{\alpha}\,t^2}{(\underline {\alpha}^2-C^2)^2\,
t^2}.\end{eqnarray*}
Finally, observe that $\frac {G(2u)-G(u)}{G(2u)-G(0)}$ tends to $\frac 
12$ as $u$ goes to zero. 
Then, by choosing $C_0$ small enough that $(\ref{C})$ holds, 
\[\frac {4C^2\bar{\alpha}\,}{(\underline {\alpha}^2-C^2)^2\,}\leq\frac 
14-p_0,\]
and that, for $u\leq C_0$, 
\[\frac {G(2u)-G(u)}{G(2u)-G(0)}\geq\frac 12-p_0,\]
the assertion is proved. 
\end{proof}

\begin{lemma}\label{tv}
\item[(i)]Let $E$ be a complete, separable metric 
space, and for $\mu_1,\mu_2\in {\cal P}(E)$ define 
\[\Vert\mu_1-\mu_2\Vert_{TV}:=\sup_{A\in {\cal B}(E)}|\mu_1(A)-\mu_
2(A)|.\]
Then there exist $\nu_0,\nu_1,\nu_2\in {\cal P}(E)$ such that 
\begin{equation}\mu_1=(1-\rho )\nu_0+\rho\nu_1,\qquad\mu_2=(1-\rho )\nu_0+\rho\nu_
2,\label{decomp}\end{equation}
where 
\[\rho =\Vert\mu_1-\mu_2\Vert_{TV}.\]
\item[(ii)]Let $E_1$ and $E_2$ be complete separable metric 
spaces and $P$ be a transition function from $E_1$ to $E_2$, and let
\[P\mu (dy)=\int_{E_1}P(x,dy)\mu (dx),\quad\mu\in {\cal P}(E_1).\]
Then, for $\mu_1,\mu_2\in {\cal P}(E_1)$ and $\nu_1,\nu_2$ as in (a), 
\[\Vert P\mu_1-P\mu_2\Vert_{TV}=\Vert\mu_1-\mu_2\Vert_{TV}\Vert P
\nu_1-P\nu_2\Vert_{TV}.\]
\end{lemma}

\begin{proof}
\noindent {\bf (i)} Let 
\[l_i:=\frac {d\mu_i}{d(\mu_1+\mu_2)},\quad i=1,2,\]
 
\begin{eqnarray*}
\nu_0(A)&:=&\frac 1{\int_E\big(l_1(x)\wedge l_2(x)\big)\,(\mu_1+\mu_
2)(dx)}\int_A\big(l_1(x)\wedge l_2(x)\big)\,(\mu_1+\mu_2)(dx),\\
\nu_i(A)&:=&\frac 1{1-\int_E\big(l_1(x)\wedge l_2(x)\big)\,(\mu_1
+\mu_2)(dx)}\int_A\big(l_i(x)-l_1(x)\wedge l_2(x)\big)\,(\mu_1+\mu_
2)(dx),\quad i=1,2,\end{eqnarray*}
and
\[\rho :=1-\int_E\big(l_1(x)\wedge l_2(x)\big)\,(\mu_1+\mu_2)(dx)
.\]
Then $(\ref{decomp})$ holds. In addition 
\[\Vert\mu_1-\mu_2\Vert_{TV}=\rho\Vert\nu_1-\nu_2\Vert_{TV}=\rho 
,\]
because $\nu_1$ and $\nu_2$ are mutually singular.

\noindent{\bf (ii)} By (i), 
\[\Vert P\mu_1-P\mu_2\Vert_{TV}=\Vert (1-\rho )P\nu_0+\rho P\nu_1
-(1-\rho )P\nu_0-\rho P\nu_1\Vert_{TV}=\rho\Vert P\nu_1-P\nu_2\Vert_{
TV}.\]
\end{proof}

\setcounter{equation}{0}

\section{Notation}\label{notation}

$\langle\cdot ,\cdot\rangle$ denotes the scalar product of two vectors. 

\noin For any matrix $M$ (or vector $v)$, $M^T$ $(v^T$) denotes 
its transpose. 

\noin$\mbox{\rm tr}M$ denotes the trace of a matrix. 

\noin For vectors $v^1,\,v^2,\,...,v^k\in\R^h$, $C(v^i,i=1,...,k)$  
denotes the closed convex cone generated by $v^1,\,v^2,\,...,v^k$.

$\I_E$ is the indicator function of a set $E$. 

\noin For $E\subseteq\R^h$ $d(x,E)$ is the distance of a point $x$ 
from $E$. 

\noin$B_r(x)\subseteq\R^h$ is a ball of radius $r$ centered at $x$.

For $f:\R^h\rightarrow\R^m$ with first order partial 
derivatives, $Df$ denotes the Jacobian matrix of $f$.

\noin For $f:\R^h\rightarrow\R$ with second order partial derivatives 
$D^2f$ denotes the Hessian matrix. 

$|\cdot |$ denotes indifferently the absolute value of a 
number, the norm of a vector or of a matrix, while 
$\|\cdot\|$ denotes the supremum norm of a bounded, reaal 
valued function. 

For an open set $E\subseteq\R^h$, ${\cal C}^i(E)$ denotes the set of 
real valued functions defined on $E$ with continuous 
partial derivatives up to the order $i$. For $E$ closed,   
${\cal C}^i(E)$ denotes the set of 
real valued functions defined on an open neighborhood of $E$ 
that admit continuous partial derivatives up to the order $i$. 
${\cal C}^i_b(E)$ denotes the subset of functions of ${\cal C}^i(
E)$ that 
are bounded with all their derivatives. 

For a complete, separable, metric space $E$, ${\cal D}_E[0,\infty 
)$ 
is the space of $E$ valued functions on $[0,\infty )$ 
that are right continuous with left hand limits, and 
${\cal C}_E[0,\infty )$ is the space of continuous functions. 

${\cal L}(\xi )$ denotes the law (distribution) of a random 
variable $\xi$.

\noin$\|\cdot\|_{TV}$ denotes the total variation norm of a 
finite, signed measure. 

Throughout the paper $c$ and $C$ denote 
positive constants depending only on the data of the 
problem. When necessary, they are indexed $c_0$, $c_1$, 
..., $C_0,$ $C_1$, ... and the dependence on the data or other 
parameters is explicitely pointed out.

\bibliography{martprob}

\end{document}